\documentclass[10pt]{article}
\usepackage{amsmath}
\usepackage{amsfonts}
\usepackage{amssymb}
\usepackage{latexsym}
\usepackage{amsbsy}
\newtheorem{definition}{Definition $\!\!$}[section]
\newtheorem{proposition}[definition]{Proposition $\!\!$}
\newtheorem{lemma}[definition]{Lemma $\!\!$}

\newtheorem{theorem}[definition]{Theorem $\!\!$}
\newtheorem{example}[definition]{\it Example $\!\!$}

\newcommand{\gm}{\gamma}
\newcommand{\om}{\omega}
\newcommand{\Om}{\Omega}
\newcommand{\Gm}{\Gamma}

\newcommand{\ot}{\otimes}
\newcommand{\sw}[1]{^{(#1)}}
\newcommand{\bpf}{{\it Proof.~~}}
\newcommand{\epf}{\hfill\mbox{$\square$}\vspace*{3mm}}
\begin{document}
{\noindent\LARGE\bf
Cup Products in Hopf-Cyclic Cohomology}\\
\medskip

\noindent{\large\bf
 Masoud Khalkhali$^a$,
 Bahram Rangipour$^b$,\\
 \medskip\medskip

\footnotesize\bf\begin{itemize}
\vspace*{-4mm}\item[$^a$]
Department of Mathematics, University of Western Ontario, London ON Canada;
E-mail: masoud@uwo.ca\\

\vspace*{-4mm}\item[$^b$]
Mathematics and Statistics,
University of Victoria,
Victoria B.C.
Canada;
E-mail : bahram@uvic.ca \end{itemize}
}
\bigskip

\ $\!\!$\hrulefill\

\medskip

\noindent{\sf Abstract.  We construct cup products of two different kinds for
Hopf-cyclic cohomology. When the Hopf algebra reduces to the ground
field our first cup product reduces to Connes' cup product in ordinary
cyclic cohomology. The second cup product generalizes
Connes-Moscovici's characteristic map for actions of Hopf algebras on
algebras.}{\small
}\medskip\ ~\\

\noindent{\sf R\'esum\'e.} {\small {\bf  Cup-produits dans la cohomologie
Hopf-cyclique}
Nous construisons deux types de cup-produits pour la
cohomologie Hopf-cyclique.
Lorsque l'alg\'ebre de Hopf se r\'eduit au corps de base,
notre premier cup-produit se r\'eduit au cup-produit de Connes
en cohomologie cyclique ordinaire.
Le deuxi\`eme cup-produit g\'en\'eralise l'application caract\'eristique
de Connes-Moscovici pour l'action des alg\`ebres de Hopf  sur les alg\`ebres.

}
\ $\!\!$\hrulefill\
\bigskip


 \section{Introduction}

 Let $H$ be a Hopf algebra, A  an
$H$-module algebra, B an $H$-comodule algebra and $C$ an
$H$-module coalgebra. Let $M$ be an stable anti-Yetter-Drinfeld
module over $H$. We denote the Hopf-cyclic cohomology of $A$, $B$,
and $C$ with coefficients in  $M$ by $HC^n_H(A, M)$,  $HC^{n,
H}(B, M)$, and   $HC^n_H(C, M)$ respectively.  Our {\it cup
product of the first kind} is a natural map
$$ HC^p_H(A, M)\ot HC^{q, H}(B, M) \longrightarrow HC^{p+q}(A\rtimes_H B),$$
with values in   the ordinary cyclic cohomology of the {\it twisted tensor
product} algebra
$A\rtimes_H B$.
For our {\it cup product of the second kind} we assume that the
coalgebra $C$ {\it acts} on the algebra $A$ in a suitable sense defined in Section 5.
Then we construct a natural map
$$HC^p_H(C, M) \otimes HC^q_H (A, M)\longrightarrow HC^{p+q}(A),$$
for all non-negative integers $p$ and $q$. Our method of construction in
both cases relies on  generalizing Connes' technique in \cite{c1}: we
realize Hopf-cyclic cocycles as {\it characters} of {\it $(H, M)$-twisted cycles} and identify those cycles
that represent a trivial cyclic cocycle.

 \section{Preliminaries}
We denote  the comultiplication,
counit, and antipode of a Hopf algebra    by $\Delta$, $\varepsilon$, and $S$ respectively.
 The letter $H$ will always
denote  a Hopf algebra over a ground field $k$ of characteristic zero.
We use the Sweedler notation, with summation sign suppressed,   to denote the
comultiplication by
   $\Delta(h)=h^{(1)}\otimes h^{(2)}$ and its higher iterations.
If $M$ is a left $H$-comodule we write $ \rho(m)=m^{(-1)} \otimes
m^{(0)},$ where summation is understood, to denote its coaction
$\rho: M \rightarrow H \otimes M$. Similarly if $M$ is a right
$H$-comodule, we write $\rho (m)= m^{(0)} \otimes m^{(1)}$ to
denote its coaction $\rho$. A left-left {\it stable
anti-Yetter-Drinfeld  (SAYD) module} is a left $H$-module and left
$H$-comodule $M$ such that
$\rho(hm)=h^{(1)}m^{(-1)}S^{-1}(h^{(3)})\ot h^{(2)}m^{(0)}$ and
$m^{(-1)}m^{(0)}=m$ for all $h\in H$ and $m\in M$.  These modules
were introduced in \cite{hkrs1, hkrs2} as the most general
coefficients systems that one can introduce into Hopf
cyclic-cohomology. One dimensional SAYD modules correspond to
Connes-Moscovici's modular pair in involution $(\delta, \sigma)$
on $H$ \cite{cm1, cm2, cm3} and will be denoted by
$M=^{\sigma}\!\!\!k_{\delta}$.

 An algebra $A$ is called a left $H$-\textit{module algebra}
 if it  is a left $H$-module and its  multiplication  and unit maps   are   morphisms of  $H$-modules.
 An algebra $B$ is called a left
 $H$-{\it comodule algebra}, if
$B$ is a left $H$-comodule and its multiplication and  unit maps  are $H$-comodule maps.
 A left   $H$-{\it module coalgebra} is a coalgebra $C$ which is
 a left $H$-module such that its comultiplication  and  counit
 maps  are
 $H$-linear.

We call the above three types of symmetries, {\it symmetries of
type A, B, and C}, respectively. For each type there is an
associated Hopf-cyclic cohomology theory with coefficients in an
SAYD $H$-module $M$ \cite{hkrs1, hkrs2}. We denote these theories
by $HC^n_H(A, M)$, $HC^{n, H}(B, M)$, and $HC^n_H(C, M)$
respectively. Connes-Moscovici's theory for Hopf algebras
\cite{cm1, cm2} correspond to the case $C=H$ with multiplication
action, and $M=^{\sigma}\!\!\!k_{\delta}$. The dual theory of
\cite{kr1} for Hopf algebras correspond to $B=H$ with
comultiplication coaction and $M=^{\sigma}\!\!\!k_{\delta}$. The
$H$-equivariant cyclic cohomology of \cite{ak} correspond to
$HC^n_H(A, M)$ with $M=H$ and conjugation action, and the twisted
cyclic cohomology with respect to an automorphism correspond to
$HC^n_H(A, M)$ with $H=k[x, x^{-1}]$ and
$M=^{1}\!\!\!k_{\varepsilon}$. $HC^0_H (A,
^{\sigma}\!\!\!k_{\delta})$ is the space of $\delta$-invariant
$\sigma$-traces on $A$ in the sense of \cite{cm2, cm3}.
\section{Covariant differential calculi}
By a differential graded (DG) left $H$-module algebra  we mean a  graded left $H$-module   algebra
$\Omega= \oplus_{i\geq 0} \,\Omega^i$ endowed with  a graded derivation
$d$
of degree 1
such that $d^2=0$ and $d$ is $H$-linear. Let $A$ be a left $H$-module algebra. Its universal
differential calculus
\cite{c1} $(\Omega A, d)$ is a DG $H$-module algebra with the $H$-action
$$h \cdot (a_0da_1\cdots da_n):=h^{(1)}(a_0)dh^{(2)}(a_1)\cdots dh^{(n+1)}(a_n).$$

By a DG $H$-module coalgebra we mean a  graded $H$-module
coalgebra $\Theta =\oplus_{i\geq 0} \Theta_i$ endowed with a
graded coderivation $d$ of degree -1
 such that  $d^2=0$ and $d$ is $H$-colinear. Let
$C$ be a left $H$-module coalgebra. Its universal differential calculus $(\Omega^c C,
d)$,
defined in\cite{k},
is a left DG $H$-module coalgebra under the $H$-action
$$h \cdot (c_0\otimes c_1 \cdots \otimes c_n)=h^{(1)}(c_0)\otimes h^{(2)}(c_1)\otimes
\cdots \otimes h^{(n+1)}(c_n).$$

By a  DG left $H$-comodule algebra we mean a DG algebra
$(\Gamma, d) $ such that $\Gamma$ is a graded left $H$-comodule
algebra and the derivation $d$ is an $H$-comodule map. If $B$ is a
left $H$-comodule algebra,  its universal calculus $(\Omega B,
d)$ is a DG left $H$-comodule algebra with the  left
$H$-coaction defined by
$$b_0db_1\cdots db_n\mapsto  (b_0\sw{-1}\cdots b_n\sw{-1})\ot b_0^{(0)}db_1^{(0)}\cdots db_n^{(0)}.$$

  \begin{definition} Let $\Omega$ be a DG left H-module algebra and $M$ be a
  left-left  SAYD module. By a closed graded $(H, M)$-trace of degree $n$  on $\Omega$ we mean a linear map
  $\int : M\otimes \Omega^n \to k$ such that for all $h\in H$, $m\in M$, and $
  \omega, \omega_1,  \omega_2$ in $\Omega$ of appropriate degrees, we
  have:
  \begin{eqnarray*}
\int h\sw{1}m\ot h\sw{2}\omega &=&\varepsilon (h)\int m\otimes
\omega,  \quad  \int m\otimes d\omega =0,\\
 \int m\otimes \omega_1 \omega_2 &=&(-1)^{deg (\omega_1) deg (\omega_2)}\int m^{(0)} \otimes
  S^{-1}(m\sw{-1})(\omega_2)\omega_1.
\end{eqnarray*}
\end{definition}

\begin{lemma} Let   $A$ be a left $H$-module algebra and
   $\rho: A \to \Omega^0$ be an $H$-linear algebra homomorphism. Then the cochain
$\varphi : M\otimes A^{\otimes (n+1)}\to k,$
$$\varphi (m, a_0, \cdots, a_n)=\int m\otimes \rho (a_0)d\rho
(a_1)\cdots d\rho (a_n),$$
is a cyclic cocycle in $Z^n_H(A, M).$ The map $\int \mapsto \varphi$ (for $\rho =id$) defines
a 1-1 correspondence between closed graded $(H, M)$-traces on $\Omega A$ and $Z^n_H(A, M)$.
\end{lemma}

\begin{definition}
Let  $\Gamma$ be a DG left  $H$-comodule algebra and $M$ be a left-left  SAYD module.
By a closed graded    $(H, M)$-trace of degree $n$ on $\Gamma$ we mean a linear map
$\int :  \Gamma^n \to M$ such that
$\int$ is $H$-colinear, and for all $m$, $\gamma, \gamma_1, \gamma_2$ of appropriate degrees
  $\int m\otimes d\gamma=0,$ and
  $$\int  \gamma_1 \gamma_2 =(-1)^{deg (\gamma_1) deg (\gamma_2)}
 \gamma_2\sw{-1}\cdot \int \gamma_2\sw{0} \gamma_1.$$
\end{definition}

\begin{lemma} Let
$B$ be a left
$H$-comodule algebra and
$\rho: B \to \Gamma^0$  an $H$-colinear algebra homomorphism.
 Then the cochain
$\varphi : B^{\otimes (n+1)}\to M,$
$$\varphi (b_0, \cdots, b_n)=\int \rho (b_0)d\rho
(b_1)\cdots d\rho (b_n),$$ is a cyclic cocycle in $Z^{n, H}(B,
M).$ The map $\int \mapsto \varphi$ (for $\rho =id$) defines a 1-1 correspondence
between closed graded $(H, M)$-traces on $\Omega A$ as in
Definition 3.3 and $Z^{n, H}(B, M)$.
\end{lemma}

\begin{definition}
Let  $\Theta $ be a DG left $H$-module coalgebra and $M$ be a left-left SAYD $H$-module. By an $n$-
dimensional closed graded  $(H, M)$-cotrace on $\Theta $ we mean  an element
$x =\sum_i m_i\ot \theta_i  \in M\ot_H\Theta_n$
such that
  $(1\ot_H d)x=0,$ and
  $$\sum_i m_i^{(0)}\otimes \theta_i^{(2)}\otimes m_i^{(-1)}\theta_i^{(1)}=\sum_i (-1)^{deg(\theta\sw{1})deg(\theta\sw{2})} m_i \otimes
  \theta_i^{(1)}\otimes \theta_i^{(2)}.$$
\end{definition}
\begin{lemma} Let $C$ be a left $H$-module coalgebra  and $\rho:
\Theta_0 \to C$ be an $H$-linear  coalgebra map. Let $\widetilde{\rho}: M\ot_H\Theta_n \to M\ot_H C^{\otimes (n+1)}$
be the natural co-extension of $\rho$. Then $\int x: =\widetilde{\rho} (x)$ is a
cyclic cocycle in $Z^n_H(C, M)$. The map $x \mapsto \int x$ (for $\rho =id$) defines a
1-1 correspondence between closed graded $(H, M)$-cotraces of degree $n$
on $\Omega^c C$ and  $Z^n_H(C, M)$.
\end{lemma}

Extending the terminology of \cite{c1}, we call the data $(\Omega, d, \int, H, M, \rho)$ an
{\it $\Omega$-cycle} over the algebra $A$ and the corresponding Hopf-cyclic cocycle its {\it character}. Similarly
for $\Gamma$-cycles and $\Theta$-cycles. They correspond to symmetries
of type $A$, $B$ and $C$, respectively.

\subsection{ Vanishing $\Omega$, $\Gamma$ and $\Theta$-cycles}
Let $B$ be a left $H$-comodule algebra and  $u\in B$ be an invertible
$H$-coinvariant element in the sense that
 $\rho (u)=1\otimes u$.
The inner automorphism
$Ad_u: B \to B$, $Ad_u(b)=ubu^{-1}$ is $H$-colinear and hence induces a
map $Ad_u^*: HC^{n, H}(B, M) \to HC^{n, H}(B, M)$
for all $n$.
\begin{proposition} We have $Ad_u^* =id$.
\end{proposition}
\bpf
One checks that the homotopy operator $\kappa$ defined by
$$\kappa f(b_0, \cdots, b_{n-1})=
\sum_{i=0}^{n} (-1)^if(b_0u^{-1}, \cdots,ub_iu^{-1}, u, b_{i+1}, \cdots, b_{n-1})$$
is $H$-colinear and is a contracting homotopy for $Ad_u^* -id$ on the  Hopf-Hochschild complex. The result
now follows by applying Connes' long exact sequence.
\epf

The algebra of $n\times n$ matrices over $B$, $M_n(B)$, is a left
$H$-comodule algebra in a natural way and the map
 $i: B \to M_n(B)$, $ b\mapsto b\otimes e_{11},$  is
$H$-colinear and hence induces a map
$i^*: HC^{p, H}(M_n(B), M) \to HC^{p, H}(B, M).$
We define a map
$Tr: C^{p, H}(B, M) \to C^{p, H}(M_n(B), M)$
\begin{eqnarray*}
(Tr \varphi)(b_0\otimes m_0,\cdots, b_p\otimes m_p)=tr (m_0\cdots
m_p)\varphi (b_0, \cdots, b_p).
\end{eqnarray*}
The relation $i^*\circ Tr=id$ is easily verified. Although we
won't need it for the construction of  cup products in this paper, we pause
to mention that  we have now all the tools to prove a  Morita invariance theorem for
Hopf-cyclic cohomology theory of any type. For example we have:
\begin{proposition} (Morita invariance)
Let $B$ be a unital left $H$-comodule algebra and $M$ be an SAYD module. Then $i^*$  induces an
isomorphism on Hopf-Hochschild and Hopf-cyclic cohomology of $B$ with coefficients in $M$.
\end{proposition}
The following lemma is an adaptation of a lemma of Connes \cite{c1} to
our context:
\begin{lemma}
Let $f: B \to B$ be an $H$-colinear
 algebra homomorphism and $X$ an invertible $H$-coinvariant
 element of $M_2(B)$ such that
 $$
  X \left[ \begin{matrix} b& 0\\ 0 & f(b) \end{matrix}\right]
  X^{-1}=\left[ \begin{matrix} 0& 0\\0& f(b)\end{matrix} \right]
 $$
 for all $b\in B$. Then for any SAYD module $M$, $HC^{n, H}(B, M)=0$ for all $n$.
\end{lemma}

\begin{definition} We say that a  $\Gamma$-cycle
is vanishing if $\Gamma^0$ satisfies the condition of the above lemma.
\end{definition}

\begin{lemma} Let $\varphi : B^{\otimes (n+1)}\to M$ be an $H$-colinear map. Then $\varphi$ is a coboundary
if and only if $\varphi$ is the character of a vanishing $\Gamma$-cycle.
\end{lemma}

 Let $C$ be a left $H$-module coalgebra and $\chi\in Hom_H(C,k)$ be an
 $H$-linear
 convolution invertible functional on $C$.
  The {\it co-inner} automorphism $Ad^c_{\chi}:C\rightarrow C$ defined by
  $Ad_{\chi}(c)=\chi(c\sw{1})c\sw{2}\chi^{-1}(c\sw{3})$ is $H$-linear and hence for
   any $n\ge 0$ induces a map
  $$ Ad^{ \ast}_{\chi}: HC^n_H(C,M)\rightarrow HC_H^n(C,M).$$
  \begin{proposition}
 We have $Ad^{ \ast}_{\chi}=id$.
 \end{proposition}
 \begin{lemma}\label{cob}
 Let  $f: C\rightarrow C$  be an $H$-linear coalgebra
 homomorphism and $\chi$  a convolution invertible $H$-linear functional  on the coalgebra
 $M_2(C)$ such that
  $$\chi (\left[\begin{matrix}c\sw{1}& 0\\ 0 & f(c\sw{1})\end{matrix}\right])\left[ \begin{matrix}c\sw{2} &0 \\ 0& f(c\sw{2})\end{matrix}\right]\chi^{-1}(\left[\begin{matrix}c\sw{3} & 0\\ 0 & f(c\sw{3})\end{matrix}\right])=\left[ \begin{matrix}0 &0\\ 0& f(c)\end{matrix}\right],$$ for all $c\in C$.
  Then for any SAYD module $M$ and for any $n\ge 0$, $HC^n_H(C,M)=0.$
 \end{lemma}
 \begin{definition}
 We say that a $\Theta$-cycle  $(\Theta,d,x,H,M)$ is vanishing if $\Theta_0$ satisfies
 the condition of Lemma \ref{cob}.
 \end{definition}

\begin{lemma} Let $y\in Z^n_H(C, M)$. Then $y$ is a coboundary
if and only if $y$ is the character of a vanishing $\Theta$-cycle.
\end{lemma}
The notion of vanishing $\Gamma$-cycle is defined along parallel lines
with analogue of Lemma 3.11 and 3. 15 proved in a similar way.

\section{Cup products of the first kind}
 In this section $A$ is a left  $H$-module algebra
 and $B$ is a left $H$-comodule algebra. Let $(\Omega, d)$ be  a DG left $H$-module algebra and
 $(\Gamma, d)$  a DG left $H$-comodule algebra.
We define
a DG algebra $\Omega \rtimes_H \Gamma$ as follows. As a vector space, it is the graded tensor product
$\Omega \otimes \Gamma$. Its multiplication and derivation are defined
by
$ (\om_1 \ot \gm_1)(\om_2 \ot \gm_2):= (-1)^{deg(\om_2)deg(\gm_1)}\om_1
\gm_1\sw{-1}(\om_2)\ot \gm_1\sw{0}\gm_2$ and
$d(\om \ot \gm):= d\om \ot\gm+(-1)^{deg(\om)}\om \ot d\gm.$  Let $M$ be a left-left SAYD module and
$\int: \Gamma \to M,  \int': M\otimes \Omega \to k,$
be closed graded $(M, H)$-traces of dimensions $p$ and $q$ on $\Gamma$
and $\Omega$, respectively.  Define a linear map
$$\int'': \Om \rtimes_H \Gm \to k, \qquad \int''(\om \ot\gm): =
\int' (\int \gm)\ot \om.$$
\begin{lemma} $\int''$ is a closed graded trace of degree $p+q$ on
 $\Om \rtimes_H \Gm$. If $\int$  (resp. $\int'$) defines  a
vanishing cycle on $B$ (resp. $A$),  then $\int''$ defines
 a vanishing
cycle on $A\rtimes_H B$.
\end{lemma}

Let $\varphi \in Z^{p, H}(B, M)$, $\psi \in Z^q_H(A, M)$ be  cyclic cocycles represented by
 $\int$
on $(\Omega B, d)$, and  $\int'$ on $(\Omega A, d)$, respectively.
Consider the sequence of maps
\begin{eqnarray*}
\Omega (A\rtimes_H B) \overset{i}{\longrightarrow} \Omega(A)\rtimes_H
\Omega(B) \overset{\int''}{\longrightarrow}k.
\end{eqnarray*}
 The first
map is induced by the natural inclusion $A\rtimes_H B \to \Om(A)\rtimes_H
\Om(B).$ We define the cup product $\varphi \# \psi$  to be the cyclic cocycle represented by the
 closed graded trace $\int'' \circ i$. Using Lemma 4.1, we obtain

\begin{theorem} With $A$, $B$, $H$, and $M$ as above, we have a natural
pairing
$$ HC^p_H(A, M)\ot HC^{q, H}(B, M) \longrightarrow HC^{p+q}(A\rtimes_H B).$$
\end{theorem}

\begin{example}{\rm Let $G$ be a discrete group acting by unital
automorphisms on an algebra $A$ and $C=H=kG$.  In \cite{kr1} the
Hopf-cyclic cohomology  of $C$  is computed in terms of
group cohomology with trivial coefficients: $HC^q_{k G}(k G,
k)=\oplus_{i\geq 0} H^{q-2i} (G).$ The cohomology groups
$HC^p_{kG}(A,k)$ are easily seen to be the cohomology of the
subcomplex of {\it invariant cyclic cochains} on $A$:
$\varphi
(ga_0, ga_1, \cdots, ga_p)=\varphi (a_0, a_1, \cdots, a_p),$ for
all $g \in G$ and $a_i \in A$.  Thus we obtain a pairing, first
considered by Connes in \cite{c2}:
$$HC^p_{kG}(A, k) \otimes H^q(G) \longrightarrow HC^{p+q}(A\rtimes
G).$$}

\end{example}

\section{Cup products of the second kind}
 In this section $C$ is a left $H$-module coalgebra and $A$ is a left
 $H$-module algebra. Let $(\Omega, d)$ be a DG left $H$-module algebra  and  $(\Theta, d)$ be a DG left
 $H$-module coalgebra.  We define the  {\it convolution DG algebra}  $Hom_H(\Theta, \Omega)$ as follows.
 As a graded vector space in degree $n$ it has $\underset{i+j=n}{\bigoplus} Hom_H(\Theta_i,\Omega^j)$. One checks
 that with convolution product
$f\ast g(\theta): =(-1)^{deg(g) deg (\theta\sw{1})}f(\theta\sw{1})g(\theta\sw{2}),$
and differential $df: =[d, f]$ (graded commutator),  $Hom_H(\Theta, \Omega)$ is a DG algebra.

Let  $\int$ be a closed graded  $(M, H)$-trace  of degree $p$ on $\Omega$  and
$x$ a closed graded    $(M, H)$-cotrace of degree $q$ on $\Theta$.
We define a  functional $\int'$ on $ Hom_H(\Theta, \Omega)$ by
$\int' f:=\int (id_M\ot f)(x).$
\begin{proposition}
 $\int'$ is a closed graded trace of degree $p+q$ on
 $ Hom_H(\Theta,\Omega).$ If  $\int$ (resp.  $x$) defines a vanishing
 cycle on $A$ (resp. $C$), then $\int'$ defines a vanishing cycle on $Hom_H(C, A)$.
\end{proposition}
Now let $\varphi \in Z^p_H(C, M)$ be represented by $ x$ and $\psi \in Z^q_H(A, M)$ by $\int$.
Consider
the sequence of maps
 $$ \Omega(Hom_H(C,A))\overset{i}{\rightarrow}  Hom_H(\Omega^c C, \Omega A)
 \overset{\int'}{\rightarrow} k,$$
 where the first map is defined  using the universal property of $\Omega$. We let $\varphi \cup \psi$
 to
 be  the character of the cycle $\int' \circ i$. Using Proposition 5.1
 we obtain
\begin{proposition}
  We  have a well-defined pairing :
$$  \cup: HC_H^p(C,M)\ot HC^q_H(A,M)\rightarrow HC^{p+q}( Hom_H (C,A)).$$
\end{proposition}

We say that  $C$ {\it acts} on  $A$ if
there is a linear map
$C\ot A\rightarrow A$
such that $c(ab)=c\sw{1}(a)c\sw{2}(b)$, $c(1)=\epsilon(c)1$,  and $(hc)(a)=h(c(a)),$
  for all  $c\in C $  $a, b\in A, h\in H.$ It is clear that the first two
  conditions are equivalent to the evaluation map
  $e: A\rightarrow Hom_H (C,A)$,    $e(a)(c)=c(a)$,  being  an algebra map.
 We thus obtain a map  $e^\ast : HC^n(Hom_H(C,A))\rightarrow HC^n(A)$. Combined with the
 above  proposition we obtain:
 \begin{theorem}
Assume  $C$ acts on $A$. Then we have a natural pairing
$$
\#=e^\ast \circ \cup : HC_H^p(C,M)\ot HC^q_H(A,M)\rightarrow
HC^{p+q}(A).$$
\end{theorem}
For $p=0$ or $q=0$ this map was already constructed in \cite{hkrs2} and
shown to coincide  (for $C=H$, $q=0$, and $M=^{\sigma}\!\!\!k_{\delta}$)
with  Connes-Moscovici's characteristic map \cite{cm1, cm2}. On the other hand, for $C=H$ and $M=^{\sigma}\!\!\!k_{\delta}$
this pairing is constructed in \cite{g}, using a different method. For a
different example,  
let $x=m\ot_Hc_0\ot c_1 \in Z_H^1(C,M)$ and $\phi\in Z^1_H(A,M)$. Then $ x\#\phi$  acts on
$A^{\ot 3}$ as follows:
$$x\#\phi (a_0\ot a_1\ot a_2)=\phi(m\ot c_0\sw{1}(a_0)\ot c_0\sw{2}(a_1)c_1(a_2))- \phi(m\ot c_0(a_0) c_1\sw{1}(a_1)\ot c_1\sw{2}(a_2)).$$
{\bf Acknowledgements}\\
We wish  to thank Professor Alain Connes for  guidance and support.

\end{document}